\documentclass[11pt]{amsart}   

\usepackage{enumerate}

\newtheorem{theorem}{Theorem}[section]
\newtheorem{lemma}[theorem]{Lemma}
\newtheorem{proposition}[theorem]{Proposition}
\newtheorem{corollary}[theorem]{Corollary}
\theoremstyle{definition}
\newtheorem{definition}[theorem]{Definition}
\theoremstyle{remark}

\numberwithin{equation}{section}

\begin{document}
\noindent                                             
\begin{picture}(150,36)                               
\put(5,20){\tiny{Submitted to}}                       
\put(5,7){\textbf{Topology Proceedings}}              
\put(0,0){\framebox(140,34)}                          
\put(2,2){\framebox(136,30)}                          
\end{picture}                                         

\def\F{{\mathcal{F}}}
\def\R{{\mathbb R}}
\def\Z{{\mathbb Z}}
\def\C{{\mathbb C}}
\def\K{{\mathbb K}}
\def\N{{\mathbb N}}
\def\P{{\mathbb P}}
\def\E{{\mathbb E}}
\def\I{{\mathbb I}}
\def\conv{{\operatorname{conv}\,}}
\def\lspan{{\operatorname{span}\,}}
\def\Supp{{\operatorname{Supp}\,}}
\def\e{{\varepsilon}}
\def\S{{\mathbb S}}

\newcommand{\diam}{\mathrm{diam}}
\def\cj#1{{\overline{#1}}}
\def\qam#1{{\widehat{#1}}}
\newcommand{\ball}[2]{\mathfrak{D}_{#2}(#1)}
\newcommand{\lball}[2]{\mathfrak{D}_{#2}^L(#1)}
\newcommand{\rball}[2]{\mathfrak{D}_{#2}^R(#1)}
\newcommand{\lnbhd}[2]{{#1}_{#2}^L}
\newcommand{\rnbhd}[2]{{#1}_{#2}^R}
\newcommand{\nbhd}[2]{{#1}_{#2}}

\def\HL{\mathcal{HL}_1}
\newcommand{\abs}[1]{\left\vert#1\right\vert}
\def\norm#1{\left\Vert#1\right\Vert}
\newcommand{\floor}[1]{\lfloor#1\rfloor}

\vspace{0.5in}

\title[Quasi-metric spaces with measure]%
{Quasi-metric spaces with measure}

\author{Aleksandar Stojmirovi\'{c}}
\address{School of Mathematical and Computing Sciences and School of Biological Sciences, Victoria University of Wellington, New Zealand}
\email{aleksand@mcs.vuw.ac.nz}
\thanks{Supported by a Bright Future PhD scholarship awarded by the NZ Foundation for Research, Science and Technology jointly with the Fonterra Research Centre and by Victoria University of Wellington and University of Ottawa research funds. The author is very grateful to his PhD supervisors Prof. Vladimir Pestov and Dr. Bill Jordan for their extensive comments and suggestions.}



\begin{abstract} 
The phenomenon of concentration of measure on high dimensional structures is usually stated in terms of a metric space with a Borel measure, also called an mm-space. We extend some of the mm-space concepts to the setting of a quasi-metric space with probability measure (pq-space). Our motivation comes from biological sequence comparison: we show that many common similarity measures on biological sequences can be converted to quasi-metrics. We show that a high dimensional pq-space is very close to being an mm-space. 
\end{abstract}

\maketitle

\section{Introduction}
\begin{definition}\label{def:qmspace}
Let $X$ be a set. A mapping $q: X \times X \rightarrow \R_+$ is called a \emph{quasi-metric} if
\begin{enumerate}[(i)]
  \item for all $x,y\in X$, \ $q(x,y)=q(y,x)=0 \iff x=y$,
  \item for all $ x,y,z\in X$, \ $q(x,z) \leq q(x,y)+ q(y,z)$.
\end{enumerate}
If $q$ is also symmetric, that is, for all $x,y\in X$, $q(x,y)=q(y,x)$, then $q$ is a metric. For each quasi-metric $q$, we denote by $\cj{q}$ its \emph{conjugate} quasi-metric, where $\cj{q}(x,y) = q(y,x)$. Furthermore, we call the metric $\qam{q}$, defined for each $x,y\in X$ by $\qam{q}(x,y) = \max\{q(x,y), q(y,x)\}=\max\{q(x,y), \cj{q}(x,y)\}$, its \emph{associated} metric. The pair $(X,q)$ is called a \emph{quasi-metric} space.

Let $w$ be a (positive) real-valued function on $X$. The triple $(X,q,w)$ is called a \emph{(generalised) weighted quasi-metric space} \cite{KuVa92, Vi99} if for all $x,y\in X$ \[q(x,y)+w(x)=q(y,x)+w(y).\] 
\end{definition}

Due to assymetry, many metric space structures naturally correspond to two quasi-metric structures, which will be henceforth referred to as the left- and right- structures wherever possible.

\begin{definition}
Let $(X,q)$ be a quasi-metric space, $x\in X$, $A, B \subseteq X$ and
$\e > 0$. Denote by
\begin{itemize}
\item $\diam(A) := \sup \{q(x,y): \ x,y \in A\}$, the diameter
of set $A$;
\item $\lball{x}{\e} := \{ y \in X: \ q(x,y) < \e\}$, the left open ball of radius $\e$ centered at $x$;
\item $\rball{x}{\e} := \{ y \in X: \ q(y,x) < \e\}$, the right open ball of radius $\e$ centered at $x$;
\item $\ball{x}{\e} := \{y \in X: \ \qam{q}(x,y) < \e\}$, the associated metric open ball of radius $\e$ centered at $x$;  
\item $q(x,A) := \inf \{q(x,y): \ y \in A\}$, the left distance from $x$ to
$A$;
\item $q(A,x) := \inf \{q(y,x): \ y \in A\}$, the right distance from $x$ to $A$;
\item $q(A,x) := \inf \{\qam{q}(x,y): \ y \in A\}$, the associated metric distance from $x$ to $A$;
\item $\lnbhd{A}{\e} := \{x\in X: \ q(A,x) < \e\}$, the left $\e$-neighbourhood of $A$;
\item $\rnbhd{A}{\e} := \{x\in X: \ q(x,A) < \e\}$, the right $\e$-neighbourhood of $A$;
\item $\nbhd{A}{\e} := \{x\in X: \ \qam{q}(A,x) < \e\}$, the associated metric $\e$-neighbourhood of $A$. 
\end{itemize}
\end{definition}

Each quasi-metric $q$ naturally induces a $T_0$ topology $\mathcal{T}(q)$ whereby a set $U$ is open if for each $x\in U$ there is $\e>0$ such that $\lball{x}{\e} \subseteq U$. The topology $\mathcal{T}(\cj{q})$ can be similarly defined in by using the right balls as its base and hence a quasi-metric space $(X,q)$ can be naturally associated with a bitopological space $(X,\mathcal{T}(q), \mathcal{T}(\cj{q}))$. Topological aspects of quasi-metric spaces have been very extensively researched - the review by K{\"u}nzi \cite{Ku01} contains 589 references! Note that a $T_0$ quasi-metric is frequently called a quasi-pseudometric \cite{Ku01} while the name quasi-metric is reserved for a map $q:X\times X$ which satisfies $q(x,y)=0 \iff x=y$ instead of axiom (i) in Definition \ref{def:qmspace}, and whose associated topology is hence $T_1$.

The main objective of this paper is to generalise various concepts related to the \emph{phenomenon of concentration of measure on high-dimensional structures} \cite{MS86,Gr99,Le01}, which are usually defined in terms of metric spaces with measure, to quasi-metric spaces with measure. While many constructions from the metric case carry through to the quasi-metric case without much change, some quasi-metric results have only trivial analogs. We will show that, in a natural sense to be defined later, 
\begin{quote}  {\it A `high-dimensional' quasi-metric space is, typically, very close to being a metric space.}
\end{quote}
Before proceeding, we will examine our motivation for doing so and in doing so provide another example of a quasi-metric space which, we believe, was not observed before.

\section{Motivation: biological sequences}

Consider sets of finite sequences over a finite alphabet $\Sigma$, denoted $\Sigma^*$. Examples of such sets are the set of all English words and, most importantly for us, sets of DNA or protein sequences. DNA sequences are formed from a four letter alphabet $\Sigma = \{A,C,G,T\}$, while the protein alphabet consists of 20 amino acids.

Search of DNA and protein sequence datasets \cite{Bai00, BKILOW03} by similarity is of fundamental importance in contemporary life sciences. The most basic search, performed using software tools such as BLAST \cite{AMSZZML97}, is the range similarity search: given a query sequence, find all the closest neighbours of that point with respect to some similarity measure.

The main similarity measure used is the Smith-Waterman \cite{SW81} \emph{local similarity score}. We will endevour to produce one of many of its equivalent definitions and show that under certain conditions, which are satisfied for most common practical cases, it can be converted to a (generalised weighted) quasi-metric.

\begin{definition}\label{defn:simscore}
Let $A\subset\N$ such that $|A|=n\in\N$. Denote by $A_i$, where $i\leq n$, the $i$-th element of $A$ (under the usual order on $\N$). If $I\subseteq\{1,2,\ldots,n\}$, set $A_I=\{A_i\in A\ |\  i\in I\}$.

Denote by $g: 2^{\N}\to\R$ a \emph{gap penalty} satisfying:
\begin{enumerate}
\item $\forall A\subset\N,\quad g(A) \geq 0$, and
\item $\forall A,B\subset\N \quad A\subseteq B \implies g(A) \leq g(B)$.
\end{enumerate}

Let $\Sigma$ be a finite alphabet. For any sequence $x\in\Sigma^n$, $n\in\N$ and any set $A\subseteq\{1,2,\ldots,n\}$, let $x_I$ denote the subsequence $x_{A_1}x_{A_2}\ldots x_{A_k}$ where $\abs{A}=k$. Let $S:\Sigma\times\Sigma\to\R$ be a map and $x\in\Sigma^m$, $y\in\Sigma^n$,$m,n\in\N$. Define the \emph{local similarity score} $s:\Sigma^*\times\Sigma^*\to\R$, by
\begin{equation*}
\label{eqn:simscore}
s(x,y) = \max_{A,\overline{A},B,\overline{B}}\{T(x_A,y_B)-g(\overline{A})-g(\overline{B})\}
\end{equation*}
where $A \subseteq \{1,2,\ldots,m\}$, $B \subseteq \{1,2,\ldots,n\}$, $|A| = |B| = k$,
$\overline{A} = \{A_1,A_1+1,\ldots,A_{k}-1,A_k\}\setminus A$, $\overline{B} = \{B_1,B_1+1,\ldots,B_{k}-1,B_k\}\setminus B$ and $T(x_A,y_B)=\sum_{i=1}^k S(x_{A_i},y_{B_i})$.
\end{definition}

The above definition can be interpreted in the following way. Firstly, two contiguous subsequences $x'$ and $y'$, of $x$ and $y$ respectively, are chosen which is why the similarity score is called local. Secondly, each letter $x'$ and $y'$ is either aligned with a letter from the other subsequence or deleted. The scores for aligned letters are given by $S$ while the costs of deletions are given by the gap penalty. Gap penalty functions may depend not only on the number of gaps but on their locations: contiguous gaps often have lower cost associated with them. Hence, we construct the local similarity score as the score of the best local alignment of two sequences given the gap penalties. 

The following result allows us to convert similarity scores to quasi-metrics. 
\begin{lemma}\label{lemma:sim2qm}
Let $X$ be a set and  $s:X\times X\to\R$ a map such that
\begin{enumerate}
\item $s(x,x) \geq s(x,y)\quad \forall x,y\in X$,
\item $s(x,y) = s(x,x) \wedge  s(y,x) = s(y,y)\implies x=y\quad \forall x,y\in X$,
\item $s(x,y)+s(y,z) \leq s(x,z) + s(y,y)\quad \forall x,y,z\in X$.
\end{enumerate}
Then $q:X\times X\to\R$ where $(x,y) \mapsto s(x,x) - s(x,y)$
is a quasi-metric. Furthermore, if $s$ is symmetric, that is, $s(x,y)=s(y,x)$ for all $x,y\in X$, $q$ is a generalised weighted quasi-metric with the weight function $w: x \mapsto s(x,x)$.
\begin{proof}
Positivity of $q$ is equivalent to (1), separation of points is equivalent to (2) while the triangle inequality is equivalent to (3). If $s(x,y)=s(y,x)$ then $q(y,x)+s(x,x) = s(y,y)-s(x,y)+s(x,x) =s(x,x)-s(x,y)+s(y,y)= q(x,y)+s(y,y)$ and thus $w: x \mapsto s(x,x)$ is a generalised weight.
\end{proof}
\end{lemma}

\begin{theorem}
Suppose $S:\Sigma\times\Sigma\to\R$ satisfies conditions of the Lemma \ref{lemma:sim2qm} and 
$S(a,a)>0$ for all $a\in\Sigma$. Then so does the similarity score $s$ on $\Sigma^*$ as defined in Definition \ref{defn:simscore}.
\begin{proof}
It is easy to see that if $S$ satisfies the Lemma \ref{lemma:sim2qm} so does $T$. Since $S(a,a)>0$ for all $x\in\Sigma$, it is clear that $s(x,x)=T(x,x)$ and thus $s(x,x)\geq s(x,y)$ for all $x,y\in X$. 

If $s(x,y)=s(x,x)$ then $s(x,y)=T(x,x)$ and hence $x$ is subsequence of $y$. Similarly, if $s(y,x)=s(y,y)$, $y$ is subsequence of $x$. Thus, $s(x,y) = s(x,x) \wedge  s(y,x) = s(y,y)\implies x=y$.

To prove the third statement pick $A,B,\overline{A},\overline{B},C,D,\overline{C},\overline{D}$ such that \begin{eqnarray*}S(x,y)& = & T(x_A,y_B)-g(x,\overline{A})-h(y,\overline{B}) \quad\text{and}\\
S(y,z)& = & T(y_C,z_D)-g(y,\overline{C})-g(z,\overline{D}).
\end{eqnarray*}
Let $I$ and $J$ be the sets of indices (possibly empty) of $A$ and $B$, and $B$ and $C$ respectively, such that $B_I=C_J=B\cap C$. It is clear that $|I|=|J|$. Denote by $K$ and $L$ the remaining indices of $B$ and $C$ respectively, that is, the sets such that $B_{K}=B\setminus C$ and $C_{L}=C\setminus B$.  

Since $T$ is a sum over sets of indices, we have
\begin{eqnarray*}
T(x_A,y_B)&=&T(x_{A_I},y_{B_I})+T(x_{A_K}, y_{B_K}) \quad\text{and}\\
T(x_C,y_D)&=&T(y_{C_J},z_{D_J})+T(z_{C_L}, z_{D_L}).
\end{eqnarray*}
Furthermore, let $\overline{A_I}$ and $\overline{D_J}$ be sets of  gaps, that is, 
\begin{eqnarray*}
\overline{A_I}&=&\{A_{I_1},A_{I_1}+1,\ldots,A_{|I|}-1,A_{|I|}\}\setminus A_I\quad\text{and}\\ 
\overline{D_J}&=&\{D_{J_1},D_{J_1}+1,\ldots,D_{|J|}-1,D_{|J|}\}\setminus D_J.
\end{eqnarray*}
Since $I$ and $J$ are subsets of indices of $A$ and $D$ respectively, $\overline{A_I}\subseteq\overline{A}$ and $\overline{D_J}\subseteq\overline{D}$ and hence $g(\overline{A_I}) \leq g(\overline{A})$ and $g(\overline{D_J}) \leq g(\overline{D})$.

Thus, $s(x,y) + s(y,z)$
\begin{eqnarray*}
= & T(x_A,y_B)-g(\overline{A})-g(\overline{B}) + T(y_C,z_D)-g(\overline{C})-g(\overline{D})\\
\leq & \quad T(x_{A_I},y_{B_I}) + T(x_{A_K}, y_{B_K}) - g(\overline{A_I}) \\
& + \, T(y_{C_J},z_{D_J})+T(y_{C_L}, z_{D_L})- g(\overline{D_J}) \\
\leq & T(x_{A_I},z_{D_J})  - g(\overline{A_I}) - g(\overline{D_J})\\ & + T(y_{B_I},y_{B_I}) + \, T(y_{B_K}, y_{B_K}) + T(y_{C_L}, y_{C_L}). 
\end{eqnarray*}

Observing that $T(x_{A_I},z_{D_J})  - g(\overline{A_I}) - g(\overline{D_J})\leq s(y,z)$ and, since $B_I$, $B_K$ and $B_L$ are disjoint subsets of indices of $y$, $T(y_{B_I},y_{B_I}) +  T(y_{A_K}, y_{B_K}) + T(y_{C_L}, y_{D_L})  \leq T(y,y) = s(y,y)$ completes the proof.
\end{proof}
\end{theorem}

The conditions of the Lemma \ref{lemma:sim2qm} are satisfied by most of the BLOSUM \cite{HH92} similarity score matrices on the amino acid alphabet, produced in the following way. Biologicaly closely related fragments of protein sequences are clustered together in the form of multiple alignments or blocks so that each row in a block represents a different fragment. The fragments within blocks are further clustered to reduce the effect of too closely related fragments and the relative frequency of observing amino acid $i$ in the same column as amino acid $j$ is denoted $\phi_{ij}$ (this is the aggregate over all columns and over all blocks). The similarity score $S$ is given by \[S(i,j) = 2\log_2\left(\frac{\phi_{ij}}{2\psi_i\psi_j}\right)\] where $\psi_i$ is the overall frequency of amino acid $i$. Hence $S$ is symmetric and it is easy to see that the triangle inequality of the quasi-metric obtained by the transformation from the Lemma \ref{lemma:sim2qm} is equivalent to \[\phi_{ij}\phi_{jk} \leq \phi_{ik}\phi_{jj}\] for all amino acids $i,j,k$. In many cases frequencies of two different amino acids being aligned are much smaller than the frequencies of amino acids being aligned with themselves and the triangle inequality is satisfied.

BLOSUM matrices are the most frequently used score matrices for similarity search of protein sequences and, as it can be seen from above, are also symmetric so that the quasi-metric obtained is generalised weighted. The similarity measures on DNA alphabet produce a metric but the distance derived from the local similarity score on DNA sequences of different length is still asymmetric. 

Quasi-metrics were investigated quite early in the development of biological sequence comparison algorithms by Waterman, Smith and Bayer \cite{WSB76}, but their emphasis at the time was on global rather than local similarity measures. Much effort was expanded on metrics \cite{Sel74a,SWF81} which were abandoned in favour of similarity score when it was realised that any `local' distance between two sequences cannot satisfy the triangle inequality. 

Most algorithms for similarity search in datasets of biological sequences, even those heuristic like BLAST \cite{AMSZZML97}, scan the whole dataset to retrieve close neighbours of a query point. Our interest is in attempting to produce indexing schemes for similarity search \cite{VPAS02,VPAS03} so that a dataset is partitioned so that very few points need to be scanned for each search. Performance of indexing schemes depends on many factors but it was observed \cite{Pestov00} that the so called `curse of dimensionality', where many indexing schemes for high-dimensional spaces perform worse than sequential scan, can be largely explained by the concentration of measure phenomenon. The results in \cite{Pestov00} refer only to metric spaces and the aim of this study is to produce foundations for studying similar phenomena in quasi-metric spaces.

It should be noted that all datasets, biological or otherwise are finite and hence topologically discrete and zero-dimensional. However, they also carry an additional structure - the normalised counting measure. Hence, each finite quasi-metric space automatically becomes a quasi-metric space with measure. 

\section{pq-spaces}

The main object of our study is the pq-space, the quasi-metric space with Borel probability measure. As two topologies can be associated with a quasi-metric, it is appropriate to use the Borel structure generated by $\mathcal{T}(q)\cup\mathcal{T}(\cj{q})$ so that any countable union, intersection or difference of any `left'- or `right'- open sets is measurable. It is easy to see that this structure is equivalent to the Borel structure generated by $\mathcal{T}(\qam{q})$, the topology of the associated metric since $\ball{x}{\e}=\lball{x}{\e}\cap\rball{x}{\e}$.

\begin{definition}
Let $(X,q)$ be a quasi-metric space, and $\mu$ a probability measure over $\mathcal{B}$, a Borel $\sigma$-algebra of measurable sets generated by $\mathcal{T}(\cj{q})$. We call the triple $(X,q,\mu)$ a \emph{pq-space}. 
\end{definition}

The pq-space is the quasi-metric analogue of the metric space with Borel measure (mm- or pm- space depending on whether the total measure is unity) defined by Gromov and Milman \cite{GrMi83,Gr99,Gr03}. For a metric space with measure, the concentration effects are expressed in terms of concentration function. Two such functions, left- and right-, can be defined for a pq-space.

\begin{definition}
Let $(X,q,\mu)$ be a pq-space and $\mathcal{B}$ the Borel $\sigma$-algebra of $\mu$-measurable sets. The \emph{left concentration function} $\alpha^L_{(X,q,\mu)}$, also denoted $\alpha^L$, is a map $\R_+\to [0,\frac{1}{2}]$ such that $\alpha^L_{(X,q,\mu)}(0)=\frac12$ and 
\[ \alpha^L_{(X,q,\mu)}(\e) = \sup\left\{1-\mu(\lnbhd{A}{\e});\ A\in \mathcal{B}, \ \mu(A) \geq \frac{1}{2}\right\}\] for $\e>0$.

Similarly, the \emph{right concentration function} $\alpha^R_{(X,q,\mu)}$, also denoted $\alpha^R$, is a map $\R_+\to [0,\frac{1}{2}]$ such that $\alpha^R_{(X,q,\mu)}(0)=\frac12$ and
\[ \alpha^R_{(X,q,\mu)}(\e) = \sup\left\{1-\mu(\rnbhd{A}{\e});\ A \in \mathcal{B}, \ \mu(A) \geq
        \frac{1}{2}\right\}\] for $\e>0$.
\end{definition}

For an mm-space $(X,d,\mu)$, $\alpha^L$ and $\alpha^R$ coincide and in that case will be denoted $\alpha_{(X,d,\mu)}$ or just $\alpha$. It is obvious that if $\diam(X)$ is finite, then for all $\e \geq diam(X)$, $\alpha^L(\e) =\alpha^R(\e) = 0$ and it can be shown that $\alpha^L$ and $\alpha^R$ are decreasing.

\begin{lemma}\label{lem:qpconcfunc}
For any pq-space $(X,q,\mu)$, for each $\e\geq 0$,
\[ 
\max\{\alpha^L_{(X,q,\mu)}(\e), \alpha^R_{(X,q,\mu)}(\e)\} \leq \alpha_{(X,\qam{q},\mu)}(\e) \leq \alpha^L_{(X,q,\mu)}(\e)+\alpha^R_{(X,q,\mu)}(\e).
\]
\begin{proof}
Let $A \in \mathcal{B}$ such that $\mu(A) \geq \frac{1}{2}$ and $\e>0$. 
Using $A_{\e}\subseteq \lnbhd{A}{\e} \cap \rnbhd{A}{\e}$,
\begin{eqnarray*}
1-\mu(\lnbhd{A}{\e}) & \leq 1 - \mu(A_{\e}) & \leq \alpha(\e) \implies \alpha^L(\e) \leq \alpha(\e) \quad \text{and}\\
1-\mu(\rnbhd{A}{\e}) & \leq 1 - \mu(A_{\e}) & \leq \alpha(\e) \implies \alpha^R(\e) \leq\alpha(\e),
\end{eqnarray*}
and it follows that $\max\{\alpha^L(\e), \alpha^R(\e)\} \leq\alpha_{(X,\qam{q},\mu)}(\e)$.

For the second inequality, use $A_{\e}\supseteq \lnbhd{A}{\e}\cap\rnbhd{A}{\e}$, and thus
$X\setminus A_{\e} \subseteq \big(X\setminus \lnbhd{A}{\e} \big) \cup \big(X\setminus \rnbhd{A}{\e} \big)$, implying
\[1-\mu(A_{\e}) \leq \big(1-\mu(\lnbhd{A}{\e}) \big)+\big(1-\mu(\rnbhd{A}{\e}) \big)\leq \alpha^L(\e) + \alpha^R(\e). \]
\end{proof}
\end{lemma}

The \emph{phenomenon of concentration of measure on high-dimensional structures} refers to the observation that in many high dimensional metric spaces with measure, the concentration function decreases very sharply, that is, that an $\e$-neighbourhood of any not vanishingly small set, even for very small $\e$, covers (in terms of the probability measure) the whole space. Examples are numerous and come from many diverse branches of mathematics \cite{Maurey79,GrMi83,AlonMilman85,MS86,Gr99,Pe02,Tal96a}. In this paper we will take a high dimensional pq-space to be a pq-space where both $\alpha^L$ and $\alpha^R$ decrease sharply. 

\subsection{Deviation Inequalities}

\begin{definition}
Let $(X,q)$  be a quasi-metric space. A map $f: X\to \R$ is called \emph{left $K$-Lipschitz} if there exists $K \in \R_+$ such that for all $x,y \in X$
\[f(x) - f(y) \leq Kq(x,y).\]
The constant $K$ is called a \emph{Lipschitz constant}.  Similarly,
$f$ is \emph{right $K$-Lipschitz} if $f(y)-f(x) \leq Kq(x,y)$. Maps
that are both left and right $K$-Lipschitz are called $K$-Lipschitz.
\end{definition}

Left 1-Lipschitz functions were studied by Romaguera and Sanchis \cite{RoSa00} under a name of semi-Lipschitz functions and used to obtain some best approximation results. We use the above terms for consistency with the remainder of our terminology. For example, it is easy to verify that the functions measuring the left or right distances to a fixed point or a set are respectively left or right 1-Lipschitz.

\begin{definition}
Let $(X,\mathcal{B},\mu)$ be a probability space and $f$ a measurable real-valued
function on $(X,q)$. A value $m_f$ is a \emph{median} or \emph{L\'evy mean} of $f$ for $\mu$
if
\[
\mu(\{f \leq m_f\}) \geq \frac{1}{2} \ \text{and} \ \mu(\{f \geq
m_f)\} \geq \frac{1}{2}.
\]
\end{definition}

A median need not be unique but it always exists. The following lemmas are generalisations of the results for mm-spaces.

\begin{lemma}\label{lemma:meddeviation}
Let $(X,q,\mu)$ be a pq-space, with left and right concentration functions $\alpha^L$ and $\alpha^R$ respectively and $f$ a left 1-Lipschitz function on $(X,q)$ with a median $m_f$. Then for any $\e > 0$  
\begin{eqnarray*}
\mu(\{x\in X: f(x) \leq m_f - \e\}) \leq \alpha^L(\e) &\ \text{and}\\ 
\mu(\{x\in X: f(x) \geq m_f + \e\}) \leq \alpha^R(\e).
\end{eqnarray*}

Conversely, if for some non-negative functions $\alpha_0^L$ and $\alpha_0^R:\R_+\to\R$,
\begin{eqnarray*}
\mu(\{x\in X: f(x) \leq m_f - \e\}) \leq \alpha_0^L(\e) &\ \text{and}\\
\mu(\{x\in X: f(x) \geq m_f + \e\}) \leq \alpha_0^R(\e)
\end{eqnarray*}
for every left 1-Lipschitz function $f:X\to \R$ with median $m_f$ and every $\e>0$, then $\alpha^L\leq\alpha_0^L$ and $\alpha^R\leq\alpha_0^R$.

\begin{proof}
Set $A=\{x\in X: f(x) \geq m_f \}$. Take any $y\in X$ such that $f(y)\leq m_f -\e$. Then, for any $x\in A$, $q(x,y)\geq f(x)-f(y)\geq\e$ and hence $q(A,y)\geq\e$, implying $y\in X\setminus \lnbhd{A}{\e}$. Therefore, $\mu(\{x\in X: f(x) \leq m_f - \e\}) \leq 1 - \mu(\lnbhd{A}{\e})\leq\alpha^L(\e)$.

Now set $B=\{x\in X: f(x) \leq m_f \}$. Take any $y\in X$ such that $f(y)\geq m_f +\e$. Then, for any $x\in B$, $q(y,x)\geq f(y)-f(x)\geq\e$ and hence $q(y,B)\geq\e$, implying $y\in X\setminus \rnbhd{B}{\e}$. Thus, $\mu(\{x\in X: f(x) \geq m_f + \e\}) \leq 1 - \mu(\rnbhd{B}{\e})\leq\alpha^R(\e)$.

The converse is equivalent to finding for each Borel set $A\subseteq X$ such that $\mu(A)\geq \frac12$, left 1-Lipschitz functions $f$ and $g:X\to\R$ with medians $m_f$ and $m_g$ respectively, such that $1-\mu(\lnbhd{A}{\e}) \leq \mu(\{x\in X: f(x) \leq m_f - \e\})$ and $1-\mu(\rnbhd{A}{\e}) \leq \mu(\{x\in X: g(x) \geq m_g + \e\})$. 

Let $A\subseteq X$ be such a set such and set for each $y\in X$, $f(y)=-q(A,y)$ and $g(y)=q(y,A)$. It is easy to see that both $f$ and $g$ are left 1-Lipschitz and that $m_f=m_g=0$. If $y\in X\setminus \lnbhd{A}{\e}$, we have $q(A,y)\geq \e$ and thus $f(y)\leq -\e$. Similarly, if $y\in X\setminus \rnbhd{A}{\e}$, we have $q(y,A)\geq \e$ implying $g(y)\geq \e$ and the result follows.
\end{proof}
\end{lemma}

Hence, we can state the alternative definitions of $\alpha^L$ and $\alpha^R$:
\[\alpha^L(\e) = \sup \big\{\mu(\{x\in X: f(x) \leq m_f - \e\}):\ f \ \text{is left 1-Lipschitz}\big\}\]
and
\[\alpha^R(\e) = \sup \big\{\mu(\{x\in X: f(x) \geq m_f + \e\}): \ f \ \text{is right 1-Lipschitz} \big\}.\]

Similar results can be easily obtained for the right 1-Lipschitz functions by remembering that if $f$ is a right 1-Lipschitz, $-f$ is left 1-Lipschitz. It is also  straightforward to observe that the absolute value of deviation of a 1-Lipschitz function from a median thus depends on both $\alpha^L$ and $\alpha^R$. 

\begin{corollary}
For any pq-space $(X,q,\mu)$, a left 1-Lipschitz function $f$ with a median $m_f$ and $\e > 0$
\[\mu(\{\abs{f-m_f } \geq\e\}) \leq \alpha^L_{(X,q,\mu)}(\e)+\alpha^R_{(X,q,\mu)}(\e).\]
\end{corollary}

This result reduces to the well-known inequality $\mu(\{\abs{f-m_f} \geq\e\}) \leq 2\alpha(\e)$ when $q$ is a metric. Deviations between the values of a left 1-Lipschitz functions at any two points are also bound by both concentration functions. 

\begin{lemma}\label{lemma:ptdeviation}
Let $(X,q,\mu)$ be a pq-space and $f\colon X\to\R$ a left (or right) 1-Lipschitz function. Then
\[(\mu\otimes\mu)(\{(x,y)\in X\times X:f(x)-f(y)\geq\e\}) \leq \alpha^L\left(\frac{\e}{2}\right) + \alpha^R\left(\frac{\e}{2}\right). \]
\begin{proof}
\begin{eqnarray*}
&&(\mu\otimes\mu) \left(\left\{(x,y)\in X\times X:f(x)-f(y) \geq \e \right\}\right)\\
&\leq& (\mu\otimes\mu) \left(\left\{(x,y)\in X\times X:f(x)-m_f\geq\frac{\e}{2}\right\}\right)\\
& + & (\mu\otimes\mu) \left(\left\{(x,y)\in X\times X: m_f - f(y) \geq \frac{\e}{2}\right\}\right)\\
&=& \mu\left(\left\{x\in X:f(x)\geq m_f +\frac{\e}{2}\right\}\right) + \mu\left(\left\{x\in X:f(x)\leq m_f -\frac{\e}{2}\right\}\right)\\
&\leq& \alpha^L\left(\frac{\e}{2}\right) + \alpha^R\left(\frac{\e}{2}\right).
\end{eqnarray*}
\end{proof}
\end{lemma}

\subsection{L\'evy families}

\begin{definition}
A sequence of pq-spaces $\{(X_n,q_n,\mu_n)\}_{n=1}^\infty$ is called \emph{left L\'evy family} if the left concentration functions $\alpha^L_{(X_n,q_n,\mu_n)}$ converge to $0$ pointwise, that is
\[\forall\e>0, \quad \alpha^L_{(X_n,q_n,\mu_n)}(\e)\to 0 \: \text{as} \ n\to\infty.\]

Similarly, a sequence of pq-spaces $\{(X_n,q_n,\mu_n)\}_{n=1}^\infty$ is called \emph{right L\'evy family} if the right concentration functions $\alpha^R_{(X_n,q_n,\mu_n)}$ converge to $0$ pointwise, that is
\[\forall\e>0, \quad \alpha^R_{(X_n,q_n,\mu_n)}(\e)\to 0 \: \text{as} \ n\to\infty.\]

A sequence which is both left and right L\'evy family will be called a L\'evy family. Furthermore, if for some constants $C_1, C_2 >0$ one has $\alpha_n(\e) < C_1 \exp(C_2 \e^2 n)$, such sequence is called \emph{normal L\'evy family}.
\end{definition}

It is a straightforward corollary of Lemma \ref{lem:qpconcfunc} that a sequence of pq-spaces $\{(X_n,q_n,\mu_n)\}_{n=1}^\infty$ is a L\'evy family if and only if the sequence of associated mm-spaces $\{(X_n,\qam{q}_n,\mu_n)\}_{n=1}^\infty$ is a L\'evy family.

To illustrate existence of sequences of pq-spaces which are right but not left L\'evy families consider the following example.

Let $X=\{a,b\}$ with $\mu(\{a\})=\frac23$ and $\mu(\{b\})=\frac13$. Set $q_n(a,b)=1$ and $q_n(b,a)=\frac1n$ where $n\in\N_+$.

It is clear that 
\[ \alpha^L_n(\e) =
\begin{cases}
\frac12 , & \text{if} \ \e = 0\\
\frac13 , & \text{if} \ 0 < \e\leq 1\\
0, & \text{if} \ \e > 1,
\end{cases}
\quad \text{and} \quad
\alpha^R_n(\e) =
\begin{cases}
\frac12 , & \text{if} \ \e = 0\\
\frac13 , & \text{if} \ 0 < \e\leq \frac1n\\
0, & \text{if} \ \e > \frac1n .
\end{cases}
\]
Hence, $\alpha^R_n$ converges to $0$ pointwise while $\alpha^L_n$ does not. In this case $\alpha_n = \alpha^L_n$. 

\section{High dimensional pq-spaces are very close to mm-spaces}

Most of the above concepts and results are generalisations of mm-space results. However, we now develop some results which are trivial in the case of mm-spaces. The main result is that, if both left and right concentration functions drop off sharply, the \emph{asymmetry} at each pair of point is also very small and the quasi-metric is very close to a metric.

\begin{definition}
For a quasi-metric space $(X,q)$, the \emph{asymmetry} is a map $\Gamma:X\times X\to\R$ defined by $\Gamma(x,y)=\abs{q(x,y)-q(y,x)}$.
\end{definition}

Obviously, $\Gamma=0$ on a metric space. However, $\Gamma$ is also close to $0$ for high dimensional spaces, that is, those pq-spaces for which both $\alpha^L$ and $\alpha^R$ decrease sharply near zero.

\begin{theorem}\label{thm:qpclosemm}
Let $(X,q,\mu)$ be a pq-space. For any $\e>0$,
\[(\mu\otimes\mu) (\{(x,y)\in X\times X:\Gamma(x,y)\geq\e\}) \leq \alpha^L\left(\frac{\e}{2}\right) + \alpha^R\left(\frac{\e}{2}\right).\] 
\begin{proof}
Fix $a\in X$ and set for each $x\in X$, $\gamma_a(x)= q(x,a)-q(a,x)$. It is clear that $\gamma_a$ is a sum of two left 1-Lipschitz maps and therefore left 2-Lipschitz. Furthermore, 
zero is its median since there is a measure-preserving bijection $(x,y)\mapsto (y,x)$ which maps the set $\{(x,y)\in X\times X: q(x,y) > q(y,x)\}$ onto the set $\{(x,y)\in X\times X: q(x,y) < q(y,x)\}$.
By lemma \ref{lemma:meddeviation}, $\mu(\{x\in X:\abs{\gamma_{a}(x)}\geq\e\}) \leq \alpha^L\left(\frac{\e}{2}\right)+\alpha^R\left(\frac{\e}{2}\right)$. Now, using Fubini's theorem,
\begin{eqnarray*}
&&(\mu\otimes\mu)(\{(x,y)\in X\times X:\abs{q(x,y)-q(y,x)}\geq\e\})\\
&= & \int_{x\in X}\int_{y\in X}
  \I_{\{\abs{\gamma_x(y)}\geq\e\}} d\mu(y) d\mu(x)\\
& \leq &  \left(\alpha^L\left(\frac{\e}{2}\right)+ \alpha^R\left(\frac{\e}{2}\right)\right)\int_{x\in X}d\mu(x)\\
& = & \alpha^L\left(\frac{\e}{2}\right)+\alpha^R\left(\frac{\e}{2}\right).
\end{eqnarray*}
\end{proof}
\end{theorem}
 
Thus, any pq-space where both $\alpha^L$ and $\alpha^R$ (and hence , by the Lemma \ref{lem:qpconcfunc}, $\alpha$) sharply decrease are, apart from a set of very small size, very close to an mm-space.

\section{Examples}
\subsection{Hamming Cube}

\begin{definition}
Let $n\in\N$ and $\Sigma=\{0,1\}$. The collection of all binary strings of length $n$ is denoted $\Sigma^n$ and called the \emph{Hamming cube}. 
\end{definition}
\begin{definition}
The \emph{Hamming distance (metric)} for any two strings $\sigma=\sigma_1\sigma_2\ldots\sigma_n$ and  $\tau=\tau_1\tau_2\ldots\tau_n \in \Sigma^n$ is given by \[d_n(\sigma, \tau) = \abs{\{i\in\N:\sigma_i\neq\tau_i\}}.\] The \emph{normalised Hamming distance} $\rho_n$ is given by \[\rho_n(\sigma,\tau) = \frac{d(\sigma,\tau)}{n} = \frac{\abs{\{i\in\N:\sigma_i\neq\tau_i\}}}{n}.\]
\end{definition}
\begin{definition}
The \emph{normalised counting measure} $\mu_n$, of any subset $A$ of a Hamming cube $\Sigma^n$ is given by \[\mu_n(A)=\frac{\abs{A}}{2^n}.\] 
\end{definition}

It is easy to see that the above definitions indeed give a set with a metric and a measure and that $(\Sigma^n, \rho_n, \mu_n)$  is an mm-space. One may wish to consider $\Sigma^n$ as a product space with $\rho_n$ as an $\ell_1$-type sum of discrete metrics on $\{0,1\}$ and $\mu_n$ an $n$-product of $\mu_1$, where $\mu_1(\{0\})=\mu_1(\{1\})=\frac12$.

The following bounds for the concentration function have been established \cite{Tal95}:

\begin{proposition}\label{prop:hammcubeconc}
For the Hamming cube $\Sigma^n$ with the normalised Hamming distance $\rho_n$ and the normalised counting measure $\mu_n$, we have \[\alpha_{(\Sigma^n, \rho_n, \mu_n)}(\e) \leq \exp(-2\e^2n).\]
\end{proposition}
Hence the sequence $\{(\Sigma^n, \rho_n, \mu_n)\}_{i=1}^{\infty}$ is a normal L\'evy family.

\subsubsection{Law of Large Numbers}
An easy consequence of the Proposition \ref{prop:hammcubeconc} is the well-known \emph{Law of large numbers}. 
\begin{proposition}\label{prop:lln}
Let $(\epsilon)_{i\leq N}$ be an independent sequence of Bernoulli random variables ($P(\epsilon_i=1)=P(\epsilon_i<=-1)=\frac12$). Then for all $t\geq 0$
\[P\left(\abs{\sum_{i\leq N} \epsilon_i} \geq t\right) \leq 2 \exp \left(-\frac{t^2}{2N}\right).\]
Equivalently, if $B_N$ is the number of ones in the sequence $(\epsilon)_{i\leq N}$ then
\[P\left(\abs{B_N-\frac{N}{2}} \geq t\right) \leq 2 \exp \left(-\frac{2t^2}{N}\right).\]
\end{proposition}

\subsubsection{Asymmetric Hamming Cube}

We will now produce a pq-space based on the Hamming cube by replacing $\rho_n$ by a quasi-metric. The simplest way is to define $q_1:\Sigma\to\R$ by $q_1(0,1)=1$ and $q_1(1,0)=q_1(0,0)=q_1(1,1)=0$ and set $q_n(\sigma, \tau)=\frac1n\sum_{i=1}^n q_1(\sigma_i,\tau_i)$. The triple $(\Sigma^n, q_n,\mu_n)$ forms a pq-space. One immediately observes that $\{(\Sigma^n, q_n, \mu_n)\}_{i=1}^{\infty}$ forms a normal L\'evy family since the associated metric $\qam{q_n}$ is the Hamming metric $\rho_n$.

Take two strings $\sigma$ and $\tau$ and let us consider the asymmetry $\Gamma_n(\sigma,\tau)$. It is easy to see that $\Gamma_n$ takes value between $0$ and $1$, being equal to the quantity \[\frac{1}{n}\Big| \abs{\{i:\sigma_i=0\wedge\tau_i=1\}} - \abs{\{i:\sigma_i=1\wedge\tau_i=0\}} \Big|.\]

Since our asymmetric Hamming cube is a product space, we can consider for each $i \leq  n$ the value $\delta_i=q(\sigma_i,\tau_i)-q(\tau_i,\sigma_i)$ as a random variable taking values of $0$, $-1$  and $1$ with $P(\delta_i=0)=\frac12$ and $P(\delta_i=-1)=P(\delta_i=1)=\frac14$ so that  $\Gamma_n(\sigma,\tau) = \frac{1}{n}\sum_{i\leq n}\abs{\delta_i}$. Now,
\begin{eqnarray*}
(\mu_n\otimes\mu_n) (\{(\sigma,\tau)\in \Sigma^n\times \Sigma^n:\Gamma_n(\sigma,\tau)\geq\e\}) & = & P\left(\sum_{i\leq n}\frac{1}{n}\abs{\delta_i} \geq\e\right) \\
& \leq & P\left(\sum_{i\leq n}\abs{\epsilon_i}\geq n\e\right) \\
& \leq & 2 \exp \left(-\frac{n\e^2}{2}\right).
\end{eqnarray*}

This is obviously the same bound as would be obtained by application of Theorem \ref{thm:qpclosemm} and Proposition \ref{prop:hammcubeconc}. 

\subsection{Penalties}

Talagrand \cite{Tal95} obtained the exponential bounds for product spaces endowed with a non-negative `penalty' function generalising the distance between two points. Penalties form a much wider class of distances than quasi-metrics but provide ready bounds for the left- and right- concentration functions.

We will outline here just one of results from \cite{Tal95} and apply it to obtain bounds for concentration functions in product quasi-metric spaces with product measure. 

Consider a probability space $(\Omega,\Sigma,\mu)$ and the product $(\Omega^N, \mu^N)$ where the product probability $\mu^N$ will be denoted by $P$. Consider a function $f:2^{\Omega^N}\times\Omega^N \to\R_+$ which will measure the distance between a set and a point in $\Omega^N$. More specifically, given a function $h: \Omega\times\Omega\to\R_+$ such that $h(\omega,\omega)=0$ for all $\omega\in\Omega$ set \[f(A,x) = \inf\left\{\sum_{i\leq N}h(x_i,y_i); y\in A\right\}.\]

\begin{theorem}[\cite{Tal95}]
Assume that
\[\norm{h}_\infty = \sup_{x,y\in\Omega}h(x,y)\] 
is finite and set 
\[\norm{h}_2 = \left(\int\int_{\Omega^2}h^2(\omega, \omega')d\mu(\omega)
d\mu(\omega')\right)^{1/2}.\]
Then
\[
P(\{f(A,\cdot)\geq u\}) \leq \frac{1}{P(A)}
\exp\left(-\min\left(\frac{u^2}{8N\norm{h}_2^2},
\frac{u}{2\norm{h}_\infty}\right)\right).
\]
\end{theorem}

If we replace $h$ above by $q_\Omega$, a quasi-metric on $\Omega$, and endow $\Omega^N$ with the $\ell_1$-type quasi-metric $q$ so that $x,y\in\Omega^N$, $q(x,y) = \sum_{i\leq N} q_\Omega(x_i,y_i)$, we have $f(A,x)=q(x,A)$ and the following corollary is obtained.

\begin{corollary}
Suppose $\norm{q_\Omega}_\infty < \infty$. Then
\[
\alpha_{(\Omega^N,q,\mu^N)}(\e) \leq 2\exp\left(-\min\left(\frac{\e^2}{8N\norm{q_\Omega}_2^2}, 
\frac{\e}{2\norm{q_\Omega}_\infty}\right)\right). 
\] 
\end{corollary}

Note that the bound applies to $\alpha$ and hence to both $\alpha^L$ and $\alpha^R$ because the norms referred to above are symmetric.

\bibliographystyle{amsplain}

\providecommand{\bysame}{\leavevmode\hbox to3em{\hrulefill}\thinspace}

\end{document}